\documentclass[letterpaper, 10 pt, conference]{ieeeconf}
\IEEEoverridecommandlockouts 
\usepackage{graphics} 
\graphicspath{ {./images/} }
\usepackage{epsfig} 
\usepackage{mathptmx} 
\usepackage{times} 
\usepackage{comment}
\usepackage{wrapfig}
\usepackage{amsmath} 
\usepackage{amssymb,color,mathtools}  
\usepackage[ruled,vlined]{algorithm2e}
\title{\LARGE \bf
Inequality Constraints in Facility Location and Other Similar Optimization Problems: An Entropy Based Approach
}
 \author{Amber~Srivastava, Gabriel Barsi Haberfeld, Naira Hovakimyan
        and~Srinivasa~M~Salapaka
\thanks{This work was supported by NSF ECCS (NRI) 18-30639}
\thanks{The authors are with the Mechanical Science and Engineering Department and Coordinated Science Laboratory, University of Illinois at Urbana-Champaign, IL, 61801 USA. E-mail: \{asrvstv6, gbh2, nhovakim, salapaka\}@illinois.edu.}
}
\begin{document}
\maketitle
\begin{abstract}
In this paper we propose an annealing based framework to incorporate inequality constraints in optimization problems such as facility location, simultaneous facility location with path optimization, and the last mile delivery problem. These inequality constraints are used to model several application specific size and capacity limitations on the corresponding facilities, transportation paths and the service vehicles. We design our algorithms in such a way that it allows to (possibly) violate the constraints during the initial stages of the algorithm, so as to facilitate a thorough exploration of the solution space; as the algorithm proceeds, this violation (controlled through the annealing parameter) is gradually lowered till the solution converges in the feasible region of the optimization problem. We present simulations on various datasets that demonstrate the efficacy of our algorithm. 
\end{abstract}

\section{INTRODUCTION}

Optimization problems such as facility location \cite{farahani2009facility}, vehicle routing \cite{toth2002vehicle}, multiway k-cut \cite{garg1994multiway}, and travelling salesman problem \cite{baranwal2017multiple} arise in many engineering applications. For instance, clustering a given dataset into different clusters, based on a similarity measure, is a widely used tool to understand and draw preliminary conclusions about the dataset \cite{gan2007data}. Similarly, several applications such as building management, battlefield surveillance, small cell network design in 5G networks \cite{ge20165g} and last mile delivery \cite{akyildiz2002survey} pose an optimization problem that requires overlaying a network of resources over the existing sensor network and designing a single or multi-hop routes from each sensor to a pre-determined destination via the network of resources \cite{kale2011maximum}. Such optimization problems are usually NP-hard \cite{aloise2009np} even in the unconstrained setting; and their complexity is further accentuated by the combinatorially large solution space. 

In this paper we expound on the optimization problems that fall into the above category and develop a framework to incorporate several {\em inequality} constraints on the underlying decision variable. In particular, we consider the (a) Facility Location Problem (FLP) \cite{cornuejols1983uncapicitated}, (b) Facility Location with Path Optimization (FLPO) \cite{kale2011maximum, srivastava2017combined} and the (b) Last Mile Delivery Problem (LMDP) \cite{cardenas2017city}, where the underlying objectives are to (a) allocate facilities to a network of spatially scattered nodes, (b) overlay a network of facilities on an existing network of nodes and design path from each node to a given destination via the network of facilities, and (c) schedule the package delivery from a transportation hub to its final destination, respectively. In several application areas, that pose the above optimization problems, the facilities, paths and the vehicles, based on their size, endurance and design capabilities, have an inherent upper bound on the number of nodes or packages they handle. For instance, in the context of FLP the retail-site selection across a city poses the problem of determining the suitable locations of the retail stores that provides easy access (in terms of travel time or distance) to all potential customers. However, owing to the infrastructural and inventory management costs few of the retail stores can only tender to the needs of a limited number of customers thereby leading to a constrained facility location problem.
\begin{figure}
    \centering
    \includegraphics[scale=0.50]{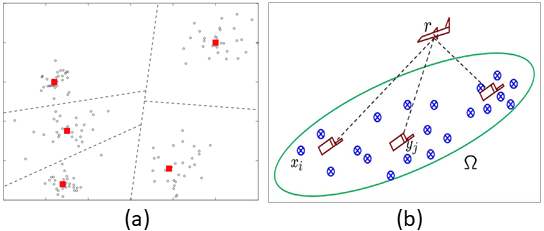}
    \caption{(a) Illustrates the FLP in $\mathbb{R}^2$ where red squares denote the $5$ allocated facilities. (b) Illustrates the FLPO problem with nodes $\{x_i\}$ that communicate with the destination at $r$ via the network of facilities $\{y_j\}$.}
    \label{fig:my_label}
\vspace{-0.5cm}
\end{figure}

Most existing algorithms address the unconstrained problems, which are complex (NP-hard) by themselves \cite{mahajan2009planar}. For instance, the FLP requires partitioning the set of $N$ nodes into $M$ clusters and allocating a facility to each of them. The number of ways in which such a partitioning can be done is combinatorially large and of the order of $^NC_M$. Also, the related cost function is non-convex with its surface riddled with multiple poor local minima. Many heuristics, such as the K-means \cite{jain2010data} and deterministic annealing \cite{rose1998deterministic} algorithms, that are used to solve the FLP, address the above issues. The FLPO problem inherits all the above complexities of FLP and comprises of the concurrent objective of designing a shortest path from each node to a given destination. In fact, in addition to combinatorially large possible partitioning of the nodes, the FLPO problem comprises of exponentially ($2^M$) large number of possible paths from each node to the destination via the network of facilities. The work done in \cite{kale2011maximum, srivastava2017combined} develop heurisitcs to solve the FLPO problem while addressing its inherent issues of non-convexity \cite{kale2011maximum} and exponentially large number of decision variables \cite{srivastava2017combined}. In the constrained setting, the equality or inequality constraints on the decision variables involved in the above optimization problems, render additional complexity to them. Several works \cite{verter2011uncapacitated,lim1999integrated, baranwal2017clustering} address the constrained optimization scenarios by developing heuristics to solve the associated integer program \cite{lim1999integrated} or adapting the existing algorithms to incorporate the constraints on an ad-hoc basis \cite{baranwal2017clustering}. 

Last Mile Delivery Problem (LMDP) concerns with the movement of the packages from a transportation hub to the final destination, both of which are usually within the same urban area. Traditionally, a vehicle is dedicated to carry the packages assigned to it all the way from the transportation hub to their respective final destinations; this results into an overall delay in the final delivery of the package as well as an expedited cost of transportation \cite{cardenas2017city}. In fact, the last mile delivery could effectively account for over $50\%$ of the package delivery cost \cite{joerss2016parcel} even though the packages are transported over comparatively very small distances.

\begin{figure}
    \centering
    \includegraphics[scale=0.4]{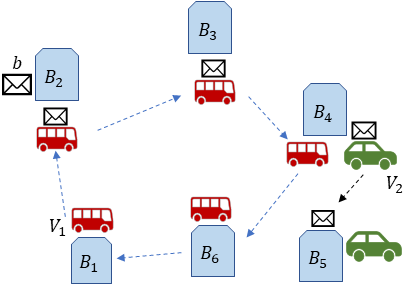}
    \caption{Illustrates the Last Mile Delivery via the network of service vehicles.}
    \label{fig: LMD}
    \vspace{-0.5cm}
\end{figure}

An alternative to the traditional last mile delivery methods is to use the {\em service} vehicles that go around the city such as the buses, metros and ride-sharing vehicles (uber and lyft) to facilitate the delivery of the packages \cite{munoz2015non, boyer2009last}. Given the time schedules and depots (the common stopping locations) for all the service vehicles, we can route the packages from the transportation hub to their final destination via these service vehicle. In particular, a package that is ready to be delivered is picked up by a service vehicle reaching the transportation hub and dropped-off either at its final destination or another suitable depot from where the package is once again picked up by another service vehicle. The underlying optimization problem is to schedule the appropriate pick-up and drop-off of each package at various depots such that they reach their respective final destination depot as quickly as possible. The Figure \ref{fig: LMD} illustrates an example scenario for the LMDP that comprises of six depots $\{B_i\}_{i=1}^6$, and two service vehicles $V_1$ and $V_2$ with routes $B_1\rightarrow B_2\rightarrow B_3\rightarrow B_4\rightarrow B_6$ and $B_4\rightarrow B_5$,  respectively. A package $b$ originates at the depot $B_2$ and is carried to the depot $B_4$ by the vehicle $V_1$, where it is again picked up by the vehicle $V_2$ and carried to its final destination depot $B_5$. For all practical purposes, the last mile delivery optimization problem involves many equality and inequality constraints that stem out of restriction on the vehicle capacity, heterogeneity in package sizes, timely delivery and handling of fragile and incompatible packages that are quite practical on the field \cite{cardenas2017city,conway2012urban}. One of our goals in this work is to model the LMDP so as to be able model and incorporate several such constraints in the associated optimization problem.

In this paper, we propose a novel Maximum Entropy Principle (MEP) \cite{jaynes1957information} based framework to address the inequality constraints in FLP, FLPO, and LMDP and present an annealing based heuristic to solve the corresponding constrained optimization problem. The MEP framework allows us to re-interpret the binary decision variables of the original optimization problem in such a way that is favourable to modelling the inequality constraints. One of our main contribution in this work is to comprehend all the inequality constraints $f_k(x)\leq F_k~\forall~k$ as an auxiliary cost function $g(x)$ such that $-$ all the inequality constraints are satisfied if and only if the auxiliary cost function attains a value below a specific upper bound (i.e., $g(x)\leq C$). We design algorithms that iteratively reduces the value of the auxiliary cost function $g(x)$ till it attains a value less than $C$ where all the inequality constraints are satisfied. In other words, our framework allows to violate the inequality constraints (i.e. $g(x)>C$) in the early stages of algorithm for better exploration of the solution space; and gradually lowers the violation till it converges to a point in the feasible region. Our proposed methodology is easily extendable to the case of equality constraints in optimization problems by re-interpreting each equality constraint as a set of two inequality constraints (i.e. $Ax=b$ is equivalent to $Ax\leq b$ and $Ax\geq b$); thus in this work we consider optimization problems with inequality constraints only. Even though we develop our framework to incorporate inequality constraints in the FLP, FLPO and LMDP, this work suggests a common framework that is applicable to many combinatorial optimization problems \cite{toth2002vehicle, jensen2011graph, baranwal2019multiway}. We demonstrate the efficacy of our proposed algorithm on randomly generated datasets in the case of FLP, FLPO and LMDP, and note that the final solution satisfies all the constraints posed upon it. 

\section{Problem Formulation and Solution}
In this section we address the optimization problems underlying the (a) FLP, (b) FLPO and (c) LMDP. We will briefly introduce the Maximum Entropy Principle (MEP) based methods \cite{rose1998deterministic,amber2018maximum} used to address the unconstrained optimization problems (a) and (b), and build upon them so as to incorporate the inequality based constraints in such problems. We then model the LMDP as a finite horizon markov decision process and develop an MEP based approach to incorporate several inequality constraints on the capacity of the service vehicles.

\subsection{Facility Location Problem}\label{sec: FLP}
The objective of the FLP is to allocate $M$ facilities $y_j\in\mathbb{R}^d,1\leq j\leq M$ for a given set of $N$ $(\gg M)$ nodes located at $x_i\in\mathbb{R}^d,1\leq i\leq N$ $-$ such that the total distance between the $N$ nodes and their closest facility gets minimized; that is, FLP aims to solve the following optimization problem 
\begin{align}\label{eq: ClusteringOrg}
\min_{\{y_j\}}\quad \sum_{i=1}^N \rho_i \big(\nu_{ij}d(x_i,y_j)\big),
\end{align}
where $\nu_{ij}=1~\text{ if }j=\arg\min_{k}d(x_i,y_k)$ else $\nu_{ij}=0$, $\rho_i$ denotes the relative importance of the $i^{\text{th}}$ node and $d(x_i,y_j)$ measures the distance between the node at $x_i$ and the facility at $y_j$. In this paper, we consider $d(x_i,y_j)$ as the squared euclidean distance cost function, i.e. $d(x_i,y_j)=\|x_i-y_j\|^2_2$. We use the MEP \cite{jaynes1957information} based DA algorithm \cite{rose1998deterministic} to address the unconstrained FLP. The DA algorithm relaxes the {\em hard} associations $\nu_{ij}\in\{0,1\}$ between the node $i$ and a facility $j$ by introducing {\em soft} associations $p_{j|i}\in[0,1]$ between the two, where $\sum_j p_{j|i} = 1$ without loss of generality. The association weights $\{p_{j|i}\}$ are designed such that they maximize the corresponding Shannon Entropy and attain a particular value $D_0$ of the cost function. In particular, we solve the following optimization problem
\begin{align}\label{eq: MEP_clustering}
\begin{split}
\max_{\{p_{j|i}\},\{y_j\}}&~H = -\sum_{i=1}^N \rho_i\sum_{j=1}^M p_{j|i}\log p_{j|i}\\
\text{subject to}&~D = \sum_{i=1}^N\rho_i\sum_{j=1}^M p_{j|i} d(x_i,y_j) = D_0.
\end{split}
\end{align}
The Lagrangian corresponding to the above optimization problem is given as
\begin{align}\label{eq: LagFLP}
F = \beta(D-D_0) - H,
\end{align}
where we refer to $F$ as the {\em free-energy} and $1/\beta$ as {\em temperature} parameter owing to its close analogies to statistical physics (where free energy is defined as enthalpy minus temperature times entropy). We minimize (local) $F$ by setting $\frac{\partial F}{\partial p_{j|i}}=0$ and $\frac{\partial F}{\partial y_j} = 0$ to obtain
\begin{align}
p_{j|i} = \frac{\exp\{-\beta d(x_i,y_j)\}}{\sum_{1\leq k \leq M}\exp\{-\beta d(x_i,y_k)\}}, ~ y_j = \frac{\sum_i \rho_i p_{j|i}x_i}{\sum_i \rho_i p_{j|i}}.
\end{align}
The constraint $D_0$ in (\ref{eq: MEP_clustering}) decides the value of the annealing parameter $\beta$. 
It is known from the sensitivity analysis \cite{jaynes2003probability} that a small $\beta(\approx 0)$ corresponds to a high value of $D_0$ and vice-versa. Also, note that at small values of $\beta(\approx 0)$, the free-energy is dominated by the Shannon Entropy H which is a convex function, and as $\beta$ increases, more and more weightage is given to the non-convex cost function $D$. The underlying idea is to determine the global minimum at $\beta = 0$ (where $F$ is convex) and track the global minimum of $F$ as $\beta$ is gradually increased. At $\beta=0$, the association weights $\{p_{j|i}\}$ are uniformly distributed and all the facilities overlap at the weighted centroid $\sum_i \rho_i x_i$ of the nodes. As $\beta$ increases, we observe no perceptible change in the location of the facilities until a critical value of $\beta = \beta_{cr1}$ is reached where the number of distinct facility locations increases. As $\beta$ increases further we observe no change in the facility locations until another critical value of $\beta=\beta_{cr2}$ is reached where the distinct facility locations once again increases. This is referred to as the {\em phase transition} phenomenon in \cite{rose1998deterministic}. Since the solution undergoes change only at these critical $\beta$'s we start the annealing process at small $\beta$ value (i.e. high $D_0$) and increase it geometrically ($\beta_{k+1}=\alpha\beta_k$, $\alpha>1$) to a large value (i.e. small $D_0$) where the number of distinct facility locations are $M$. As $\beta\rightarrow \infty$ the free-energy $F$ converges to the non-convex function $D$ and the associations become {\em hard} i.e. the association weights $p_{j|i}\rightarrow \{0,1\}$. In fact, one can explicitly compute the value of critical $\beta$'s using the necessary conditions of optimality namely, (a) $\partial F/\partial Y=0$ and (b) $d^2F/dY^2> 0$ where the phase transition occurs when the Hessian $\partial F^2/\partial Y^2$ loses rank \cite{rose1998deterministic}.

Now we consider the class of problems where the facilities $\{f_j\}_{j=1}^M$ have an inherent capacity constraint on the fraction $\{c_j\}_{j=1}^M$ of total nodes that are associated to them. Here we propose our novel MEP based framework wherein we model the inequality constraints in terms of the decision variables introduced in the MEP and re-interpret the inequality constraints as an auxiliary cost function in our optimization problem. Note that the effective fraction $p_j$ of nodes that are associated to the $j^{\text{th}}$ facility is given by $p_j=\sum_{i=1}\rho_i p_{j|i}$; thus the inequality constraints are
\begin{align}\label{eq: ineq_constraint}
p_j \leq c_j~\forall~1\leq j\leq M.
\end{align}
We consider the auxiliary cost function $\sum_j \exp\big\{\theta(p_j-c_j)\big\}$ where $\theta \gg 1$. Note that only when all the inequality constraints are satisfied the above auxiliary cost function attains a small value $\ll 1$ otherwise it attains a large value ($\gg 1$). In our MEP framework, we add the auxiliary cost function to the Lagrangian \ref{eq: LagFLP} as an equality constraint that requires $\sum_j \exp\{\theta(p_j-c_j)\}=\mu$ and design our algorithm that solves the consecutive optimization problems at gradually decreasing values of $\mu$ till it attains a value at which all the inequality constraints in (\ref{eq: ineq_constraint}) are satisfied. In the constrained FLP we seek to minimize the Lagrangian
\begin{align}
\min_{\{p_{j|i}\},\{y_j\}}\bar{F} = F + \beta'\big(\sum_{j=1}^M\exp\big\{\theta(p_j-c_j)\big\}-\mu\big),
\end{align}
where $F$ is given in (\ref{eq: LagFLP}) and $\beta'$ is a Lagrange multiplier. We minimize (local) the free-energy $\bar{F}$ by setting $\frac{\partial \bar{F}}{\partial p_{j|i}}=0$ and $\frac{\partial \bar{F}}{\partial y_j}=0$ to obtain
\begin{align}\label{eq: weight_loc_constraint}
\begin{split}
p_{j|i} &= \frac{e^{-\beta d(x_i,y_j)-\beta'\theta e^{\theta(p_j-c_j)}}}{\sum_{k=1}^M\exp^{-\beta d(x_i,y_k)-\beta'\theta e^{\theta(p_k-c_k)}}}~
y_j = \frac{\sum_{i}\rho_i p_{j|i}x_i}{\sum_{i}\rho_i p_{j|i}}.
\end{split}
\end{align}
It is known from sensitivity analysis \cite{jaynes2003probability} that a small value of $\beta'$ corresponds to a large value of $\mu$ and vice-versa. We design our algorithm in such a way that for every $\beta$ we gradually increase $\beta'$ (equivalently, we decrease $\mu$), until $\mu$ reaches some appropriate small value where all the constraints in (\ref{eq: ineq_constraint}) are satisfied. As in the unconstrained case, here also we observe the phase transitions at $\bar{\beta}_{cr}=(\beta_{cr},\beta_{cr}')$ where the number of distinct facility locations increases; whereas no significant change is observed between two consecutive critical $\bar{\beta}_{cr}$'s. Thus, we design geometric annealing laws for $\beta$ and $\beta'$, i.e. $\beta\leftarrow\alpha\beta$ and $\beta'\leftarrow\alpha'\beta'$ where $\alpha,\alpha'>1$; thereby making our proposed algorithm computationally efficient. As a part of our ongoing research we are working on determining explicit values of $\bar{\beta}_{cr}$. Please refer to the Algorithm \ref{Alg: Algorithm1} for implementation details.
\begin{algorithm}
\KwIn{$\{x_i\}_{i=1}^N$, $\{c_j\}_{j=1}^M$,  $\beta_{1\min},\beta_{1\max},\beta_{2\min},\beta_{2\max}$, annealing rates $\alpha_1>1$,and $\alpha_2>1$.}
\textbf{Initialize: $\beta_1 = \beta_{1\min}$, $y_j = \sum_i \rho_i x_i~\forall~j$}\\
 \While{$\beta_1\leq \beta_{1\max}$}{$\beta_2 = \beta_{2\min}$\\
  \While{$\beta_2\leq\beta_{2\max}$}{
  Solve the implicit equations in (\ref{eq: weight_loc_constraint}); $\beta_2 \leftarrow \alpha_2\beta_2$
  }
  $\beta_1 \leftarrow \alpha_1\beta_1$
 }
 \caption{DA Algorithm with Inequality Constraints}\label{Alg: Algorithm1}
\end{algorithm}

\subsection{Facility Location with Path Optimization}
The Facility Location with Path Optimization (FLPO) problem \cite{amber2018maximum} involves a two-fold objective of (a) allocating facilities $\{f_j\}_{j=1}^M$ to a network of nodes $\{n_i\}_{i=1}^N$ and (b) determining single or multi-hop path from each node $n_i$ to a given destination $\delta$ via the network of allocated facilities $-$ such that the sum total of cost incurred along all the paths from the nodes to the destination gets minimized. Here the node $n_i$ is located at $x_i\in\mathbb{R}^d$ $\forall$ $1\leq i \leq N$, the destination $\delta$ is located at $z\in\mathbb{R}^d$ and the location of the facility $f_j$ is denoted by $y_j\in\mathbb{R}^d$ $\forall$ $1\leq j\leq M$. As in \cite{amber2018maximum}, a path is defined as a sequence $\gamma=(\gamma_1,\hdots,\gamma_M)$ of $M$ steps where each step $\gamma_k$ corresponds to either one of the facilities $f_j$ or to the destination $\delta$; in other words $\gamma_k\in \{f_1,f_2,\hdots,f_M,\delta\}$ $\forall$ $1\leq k\leq M$ and the path $\gamma$ from the node $n_i$ to the destination $\delta$ is illustrated as
$n_i\rightarrow \gamma_1\rightarrow\gamma_2\hdots\rightarrow\gamma_M\rightarrow\delta$. The objective of the FLPO problem is
\begin{align}\label{eq: FLPO_costFunc}
\min_{\{y_j\}}~ \sum_{i=1}^N \rho_i\big(q_{i\gamma}d(n_i,\gamma)\big)
\end{align}
where $q_{i\gamma}=1\text{ if }\gamma=\arg\min_{\gamma'\in\mathcal{G}}d(i,\gamma')$ else $q_{i\gamma}=0$, $\mathcal{G}:=\{(\gamma_1,\hdots,\gamma_M):\gamma_k\in\{f_1,\hdots f_M,\delta\}~\forall~k\}$  denotes the set of all possible paths, 
$\rho_i$ denotes the relative importance of the node $n_i$, and $d(n_i,\gamma)$ is the cost incurred from the node $n_i$ to the destination $\delta$ along the path $\gamma = (\gamma_1,\hdots,\gamma_M)$; more specifically,
$d(n_i,\gamma) = d_0(n_i,\gamma_1)+ d_1(\gamma_1,\gamma_2) + \hdots + d_M(\gamma_M,\delta)$
where $d_k(\cdot,\cdot)$ is considered to be the square euclidean cost between the steps $\gamma_k $ and $\gamma_{k+1}$. For instance if $\gamma_k = f_j$ and $\gamma_{k+1} = f_{j'}$ then $d_k(f_j,f_{j'}) = \|y_j-y_{j'}\|_2^2$. We begin with replacing the {\em hard} association $q_{i\gamma}\in\{0,1\}$ between a node $n_i$ and a path $\gamma\in\mathcal{G}$ with the {\em soft} association $p(\gamma|i)\in[0,1]$ where without loss of generality we assume $\sum_{\gamma\in\mathcal{G}}p(\gamma|i)=1$ $\forall$ $i$. We design these association weights $\{p(\gamma|i)\}$ such that they maximize the Shannon Entropy while attaining a specified value of the relaxed cost function $D_p$ as described below
\begin{align}\label{eq: MEP_flpo}
\begin{split}
\max_{\{p(\gamma|i)\},\{y_j\}}&~H_p := -\sum_{i=1}^N \rho_i\sum_{j=1}^M p_{j|i}\log p_{j|i}\\
\text{subject to}&~D_p:=\sum_{i=1}^N\rho_i\sum_{j=1}^M p(\gamma|i)d(n_i,\gamma) = D_{p0}.
\end{split}
\end{align}
The corresponding Lagrangian $F_p$ for the above optimization problem is
\begin{align}\label{eq: LagFLPO}
F_p = \beta(D_p -D_{p0}) - H_p.
\end{align}

Note that the {\em law of optimality} enables to dissociate the weight $p(\gamma|i)$ into the product of step-wise association weights $p_k(\gamma_{k+1}|\gamma_k)$ for $0\leq k \leq M-1$. More precisely, for $\gamma_0=i$,
\begin{align}
p(\gamma|\gamma_0) = p_0(\gamma_1|\gamma_0)p_1(\gamma_2|\gamma_1)\hdots p_{M-1}(\gamma_M|\gamma_{M-1}).
\end{align}
We minimize (local) the Lagrangian $F_p$ by setting $\frac{\partial F_p}{\partial p_k(\gamma_{k+1}|\gamma_k)}=0$ and $\frac{\partial F}{\partial y}=0$ to obtain the expressions for $\{p_k(\gamma_{k+1}|\gamma_k)\}$ and $\{y_j\}$.
The algorithm proposed in \cite{amber2018maximum} for the unconstrained FLPO problem demonstrates the traits similar to DA algorithm in Section \ref{sec: FLP}. In particular, we observe that as $\beta$ is increased, at certain critical $\beta$'s the algorithm undergoes {\em phase transition} where the number of distinct facility location increases; and for all other $\beta$ values there is no perceptible change in the facilities. As $\beta\rightarrow \infty$ the Lagrangian $F_p\rightarrow D_{p}$ and we obtain hard associations i.e. $p_k(\gamma_{k+1}|\gamma_k)\in\{0,1\}$. 

Now we move on to the case of constrained Facility Location and Path Optimization problem where each facility $f_j$ has a given capacity $w_j$ which upper bounds the fraction of nodes $n_i$ that avail the services of $f_j$ in any of the steps $k\in\{1,\hdots,M\}$. Note that in the MEP framework for FLPO problem the effective fraction of nodes $\{n_i\}$ that a facility $f_j$ caters to is given by the expression
\begin{align}\label{eq: measureCap}
C(f_j) &= \sum_{i=1}^N \rho_i p_0(f_j|i) + \sum_{i=1}^N \sum_{\gamma_1}\rho_i p_0(\gamma_1|i)p_1(f_j|\gamma_1)\nonumber\\
&~ +\hdots+\sum_{i=1}^N\sum_{\gamma_1,\hdots,\gamma_M}\rho_i p_0(\gamma_1|i)\hdots p_{M-1}(f_j|\gamma_{M-1}),
\end{align}
where the first term in (\ref{eq: measureCap}) measures the fraction of nodes that avail the services of the facility $f_j$ in the first step of their corresponding path to the destination. Similarly, the second term measures the fraction of nodes that $f_j$ handles in the second step of their respective paths to the destination and so forth till the last term in (\ref{eq: measureCap}) which measures the fractions of nodes that $f_j$ caters to in the last step of their respective paths the destination $\delta$. Thus the inequality constraints posed by the constrained FLPO problem are
\begin{align}\label{eq: FLPO_const}
C(f_j) \leq w_j~\forall~j\in\{1,\hdots,M\}.
\end{align}
Analogous to our method in Section \ref{sec: FLP}, we choose the auxiliary cost function as $\sum_j \exp\big\{\theta(C(f_j)-w_j)\big\}$ where $\theta \gg 1$. Note that only when all the inequality constraints are satisfied the above auxiliary cost function attains a small value $\ll 1$ otherwise it takes up a large value $\gg 1$. In the MEP framework we add the auxiliary cost function to the Lagrangian $F_p$ in (\ref{eq: LagFLPO}) as an equality constraint which requires that $\sum_{j=1}^M \exp\{\theta(C(f_j)-w_j))\}=\mu_p$ and design our algorithm that solves the consecutive optimization problems at gradually decreasing values of $\mu_p$ till it reaches a value corresponding to which all the inequality constraints in (\ref{eq: FLPO_const}) are satisfied. In other words, we seek to minimize the Lagrangian 
\begin{align}
\bar{F}_p = F_p + \beta'\big(\sum_{j=1}^M e^{\theta(C(f_j)-w_j)})-\mu_p\big),
\end{align}
where $F_p$ is given in (\ref{eq: LagFLPO}), $\beta'$ is a Lagrange parameter. We exploit the inverse correlation between $\mu_p$ and $\beta'$ in our algorithm wherein for each $\beta$ we reduce the value of $\mu_p$ by analogously increasing the annealing parameter $\beta'$. In particular, we perform the two steps (a) for every given $\beta$ value we gradually increase $\beta'$ from a small value (equivalently, large $\mu_p$) to a large value (equivalently, small $\mu_p$) and (b) gradually increase $\beta$ (decrease $D_{p0}$ in (\ref{eq: MEP_flpo})) from a small to a large value (equivalently $D_0$ goes from a large to a small value); and exploit the underlying phase transitions to allow geometric annealing laws for both $\beta$ and $\beta'$. Note that as $\beta,\beta'\rightarrow\infty$, $F_p$ converges to the original cost function $D_p$ and the association weights $p_k(\gamma_{k+1}|\gamma_k)\rightarrow\{0,1\}$. Please refer to the Algorithm \ref{Alg: Algorithm2} for implementation details.

We minimize (local) the free-energy $\bar{F}_p$ with respect to $p_k(\gamma_{k+1}|\gamma_k)$ and $Y:=\{y_j\}$ by setting $\frac{\partial \bar{F}_p}{\partial p_k(\gamma_{k+1}|\gamma_l)}=0$ and $\frac{\partial \bar{F}_p}{\partial Y}=0$ to obtain
\begin{align}\label{eq: FLPO_cap_constSoln}
p_k=e^{-\bar{d}_k}\frac{\mathlarger{\smashoperator[r]{\sum_{\substack{(\gamma_{k+2},\hdots,\gamma_M)}}}e^{\sum_{t=k+1}^M -\bar{d}_t}}}{\mathlarger{\smashoperator[r]{\sum_{\substack{(\gamma_{k+1},\hdots,\gamma_M)}}}e^{\sum_{t=k}^M-\bar{d}_t}}},\quad
y = (2\hat{A}-\hat{B})^{-1}(\hat{\bar{X}} + \hat{C}),
\end{align}
where $\bar{d}_t=\beta d_t(\gamma_t,\gamma_{t+1})+\theta\beta'e^{\theta[C(\gamma_{t+1})-w_{\gamma_{t+1}}]}$
where $p_k = p_k(\gamma_{k+1}|\gamma_k)$, $\hat{A}=I_2\otimes A$, $\hat{B} = I_2\otimes B$, $\hat{C}=I_2\otimes C$, $\hat{\bar{X}}=I_2\otimes \bar{X}$, and $A,B\in\mathbb{R}^{M\times M}$, $\bar{X},C\in\mathbb{R}^{M\times d}$ depend on the association weights $p_k(\gamma_{k+1}|\gamma_k)$, the spatial locations of the node, facilities and the destination, as illustrated in the Appendix.

\begin{algorithm}
\KwIn{$\{x_i\}_{i=1}^N$, $\{c_j\}_{j=1}^M$,  $\beta_{\min},\beta_{\max},\beta_{\min}',\beta_{\max}'$, annealing rates $\alpha_1>1$,and $\alpha_2>1$.}
\textbf{Initialize: $\beta = \beta_{\min}$, $y_j = \sum_i \rho_i x_i~\forall~j$}\\
 \While{$\beta\leq \beta_{\max}$}{
 $\beta' = \beta_{\min}'$\\
  \While{$\beta'\leq\beta_{\max}'$}{
  Solve the implicit equations in (\ref{eq: FLPO_cap_constSoln}); $\beta \leftarrow \alpha_2\beta'$
  }
  $\beta\leftarrow \alpha_1\beta$
 }
 \caption{FLPO Problem with Inequality Constraints}\label{Alg: Algorithm2}
\end{algorithm}

\subsection{Last Mile Delivery Problem}
In this section we first model the last mile delivery problem involving service vehicles as a finite horizon MDP and pose it as an optimization problem in the MEP framework to determine optimal package delivery schedules. We then extend our framework to model and incorporate capacity constraints on the service vehicles. As in the case of FLP and FLPO, we think of the inequality constraints as an auxiliary cost function that is suitably minimized to attain a value corresponding to which all the inequality constraints are satisfied. Let $\mathcal{B}=\{B_1, B_2,\hdots,B_M\}$ denote all the depots, $\mathcal{V}=\{V_1,V_2,\hdots,V_N\}$ denote the service vehicles and $\mathcal{P}=\{b_1,\hdots,b_R\}$ denote the set of all packages. For each vehicle $V_k$, $1\leq k\leq N$ the route information is given in the form of sequence of depots $B_{l_1}(V_k)\rightarrow B_{l_2}(V_k)\rightarrow \hdots\rightarrow B_{l_{q_k}}(V_k),\text{ where } B_{l_r}(V_k)\in \mathcal{B}~\forall r$
that the vehicle $V_k$ visits and the corresponding times $t(B_{l_1}(V_k))\rightarrow t(B_{l_2}(V_k))\rightarrow\hdots\rightarrow t(B_{l_{q_k}}(V_k))$ of leaving each such depot. For each package $b_j$, $1\leq j \leq R$ the origin depot $B_o(b_j)$ and the destination depot $B_d(b_j)$ are given. The objective is to schedule the pick-up and drop-off of each package at the appropriate depots by the service vehicles $-$ such that the total time taken by each package to reach its destination depot is minimized. 

\begin{figure}
    \centering
    \includegraphics[scale=0.30]{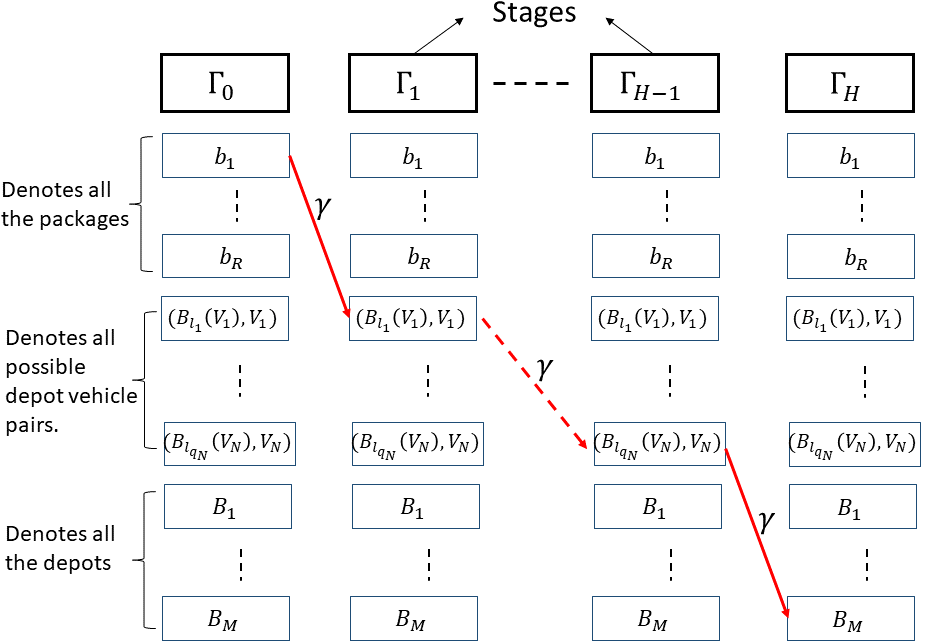}
    \caption{Stage-wise illustration of the LMDP. The corresponding finite MDP comprises of $M=|\mathcal{S}|$ states. Owing to a deterministic state transition probability $P$, any package $b_j\in\mathcal{P}$ reaches its destination depot in a maximum of $H$ stages (or steps) provided there exists atleast one feasible route for the package. Each stage $\Gamma_t,0\leq t\leq H$ comprises of all the states in $\mathcal{S}$. The state $s=(B_{l_r}(V_k),V_k)$ indicates the vehicle $V_k$ leaving the depot $B_{l_r}(V_k)$. The path $\gamma$ (in red) from $s_0=b_1$, as denoted above, indicates that the package $b_1$ is picked up by the vehicle $V_1$ from the depot $B_{l_1}(V_1)$ in the first stage $\Gamma_1$. Subsequently, it reaches the state $s_{H-1}=(B_{l_{q_N}}(V_N),V_N)$ or equivalently it reaches the depot $B_{l_{q_N}}$ from where the vehicle $V_N$ takes it to its destination depot $B_M$. Note that the cost $c(s,a)$ is computable easily from the given vehicle schedules. In case it is not possible to go from state $s_t\in\Gamma_t$ to $s_{t+1}\in\Gamma_{t+1}$, then the associated cost $c(s=s_t,a=s_{t+1})$ is assumed infinite.} 
    \label{fig:LMDP_explain}
\vspace{-0.5cm}
\end{figure}

We model this problem as a finite MDP $\mathcal{M}=\langle \mathcal{S},\mathcal{A},c,P,H\rangle$ where $\mathcal{S}=\{b_j,(B_{l_r}(V_k),V_k),B_i:1\leq j\leq R,1\leq r\leq q_k,1\leq k\leq N,1\leq i \leq M\}$ denotes the state space (see Figure \ref{fig:LMDP_explain}), $\mathcal A = \mathcal S$ denotes the action space such that an action $a\in \mathcal A$ from the current state $s\in\mathcal{S}$ takes the system to the state $s'=a\in \mathcal{S}$ while incurring a cost as defined by the function $c:\mathcal{S}\times \mathcal{A}\rightarrow \mathbb{R}$. The function $c(\cdot,\cdot)$ denotes the time taken to go from the state $s$ to the next state $s'=a$; this is explicitly calculable using the route and schedule information of the vehicles, $P:\mathcal{S}\times\mathcal{S}\times\mathcal{A}\rightarrow [0,1]$ denotes the deterministic state transition probability matrix defined as $P(s'|s,a)=1$ if $s'=a$, $0$ otherwise, $H = |\mathcal{S}|$ denotes the horizon of the finite MDP. Figure \ref{fig:LMDP_explain} provides a stage-wise illustration of the finite MDP where each stage $\Gamma_k~,0\leq k \leq H$ comprises of all the states in the state space $\mathcal{S}$. Let $\gamma\in\mathcal{G}:=\{(a_0,\hdots,a_{M-1}),~a_h\in \mathcal{A},0\leq h\leq M-1\}$ denote the sequence of actions taken from an initial states $s_0\in\mathcal{S}$; and the corresponding cost incurred be $\sum_{h=0}^{M-1} c(s_h,a_h)$ where $s_{h+1}=a_h,~0\leq h\leq M-1$. The objective in the context of LMDP is to minimize the cost function
\begin{align}\label{eq: costfuncLMD}
&\sum_{b_j\in\mathcal{P}}\rho_{b_j} \mu_{b_j\gamma}\sum_{h=0}^{M-1} c(s_h,a_h)~\text{where }s_0=b_j
\end{align}
where $\mu_{b_j\gamma}=1$ if $\gamma=\arg\min_{(a_0,\hdots,a_{M-1})}\sum_{h=0}^{M-1}c(s_h,a_h)$ else $\mu_{b_j\gamma}=0$ and $\rho_{b_j}$ is the relative importance of the package $b_j$. In the stage-wise illustration of the finite MDP in Figure \ref{fig:LMDP_explain} the above optimization problem is equivalent to designing shortest path from each $b_j\in\Gamma_0$ to their respective destination depot in $\Gamma_{H}$. For instance, $\gamma$ (indicated in red) in Figure \ref{fig:LMDP_explain} illustrates the route taken by the package $b_1$ where it goes to the state $B(l_1(V_1),V_1)\in\Gamma_1$ (or equivalently is picked up by the vehicle $V_1$ from the originating depot $B_0(b_1)=B_{l_1}(V_1)$) and is subsequently taken to its destination depot $B_M$ from the depot $B_{l_{q_N}}(V_N)$ by the service vehicle $V_N$.
We use the MEP framework to address the optimization problem in (\ref{eq: costfuncLMD}) where we relax the cost function by replacing the hard associations $\mu_{b_j\gamma}\in\{0,1\}$ with soft associations $p(\gamma|b_j)\in [0,1]$ and then use the law of optimality to re-interpret $p(\gamma|b_j)$ in terms of the stage-wise association weights $p_{h}(\gamma_{h+1}|\gamma_h)$ as
\begin{align}
p(\gamma|b_j) = \prod_{h=0}^{M-1}p_h(\gamma_{h+1}|\gamma_h)
\end{align}
\hspace{-0.15cm}where $\gamma\in\mathcal{G}, \gamma_0 = b_j, \gamma_h = a_h~ \forall h\geq 1$. The MEP poses the following optimization problem
\begin{align}\label{eq: eqn2}
\begin{split}
\max_{\{p_h(\gamma_{h+1}|\gamma_h)\}}& H_M = -\sum_{b_j\in\mathcal{P}}\sum_{\gamma\in\mathcal{G}}p(\gamma|b_j))\log p(\gamma|b_j)\\
& D_M=\sum_{b_j\in\mathcal{P}}\sum_{\gamma\in\mathcal{G}}p(\gamma|b_j)\sum_{h=0}^{M-1}c(s_h,a_h)=D_{M0}
\end{split}
\end{align}
that results into the Gibbs distribution
\begin{align}\label{eq: eqn1}
p_h(\gamma_{h+1}|\gamma_h)=e^{-\beta c_h(\gamma_{h},\gamma_{h+1})}\frac{\mathlarger{\smashoperator[r]{\sum_{\substack{(\gamma_{h+2},\hdots,\gamma_M)}}}e^{-\beta\sum_{t=h+1}^M c_t(\gamma_t,\gamma_{t+1})}}}{\mathlarger{\smashoperator[r]{\sum_{\substack{(\gamma_{h+1},\hdots,\gamma_M)}}}e^{-\beta\sum_{t=h}^M d_t(\gamma_t,\gamma_{t+1})}}},
\end{align}
where $\beta$ is the Lagrange parameter corresponding to the constraint in (\ref{eq: eqn2}). The Lagrangian $F_M$ is a convex function of the association weights $\{p_h(\gamma_{h+1}|\gamma_h)\}$ $\forall$ $\beta$, thus in our algorithm we directly set $\beta\rightarrow\infty$ to obtain $p_k(\gamma_{k+1}|\gamma_k)\in\{0,1\}$ that minimizes the cost function $D_M$. 

Now we consider the class of problems where the service vehicles have an associated capacity constraint. In particular, the service vehicles $\{V_k\}_{k=1}^N$ have an associated upper bound on the fraction $\{w_k\}_{k=1}^N$ of total packages they can carry. Note that in the MEP framework, the effective fraction of packages that a vehicle $V_k$ carries from a depot $B_{l_r}(V_k)$ is given by the following expression
\begin{small}
\begin{align}
&C\big(\bar{V}_k\big)=\sum_{b_j\in\mathcal{P}}\rho_j p_0(\bar{V}_k|b_j) +\smashoperator[r]{\sum_{{b_j\in\mathcal{P}, \gamma_1\in A}}}\rho_j p_1(\gamma_1|b_j)p_2(\bar{V}_k|\gamma_1)+\hdots\nonumber\\
&+ \smashoperator[r]{\sum_{\substack{b_j\in\mathcal{P},\gamma_1,\hdots,\gamma_{M-1}}}}\rho_jp_1(\gamma_1|b_j)\hdots p_{M-2}(\gamma_{M-1}|\gamma_{M-2})p_{M-1}(\bar{V}_k|\gamma_{M-1})
\end{align}
\end{small}
\hspace{-0.10cm}where $\bar{V}_k := \big(B_{l_r}(V_k),V_k\big)$. The first term in the above expression measures the fraction of total packages that originate at the depot $B_{l_r}(V_k)\in \mathcal{B}$ and are picked up by the vehicle $V_k$, the second term corresponds to the fraction of total packages picked up from the depot $B_{l_r}(V_k)$ by the vehicle $V_k$ where for all such packages $B_{l_r}(V_k)$ is the {\em second} depot en-route to their final destination depot. Similarly, the last term measures the fraction of total packages picked up from the depot $B_{l_r}(V_k)$ by the vehicle $V_k$ where for all such packages $B_{l_r}(V_k)$ is the $M-1^{\text{th}}$ or {\em penultimate }depot en-route to their final destination depot. Thus, the upper bound of $w_k$ on the fraction of packages carried by the service vehicle $V_k$ at any instant along its route is equivalent to the following set of inequality constraints
\begin{align}\label{eq: constLMDP}
C\big(B_{l_r}(V_k),V_k\big) \leq w_k ~\forall~1\leq r\leq q_k.
\end{align}
In other words, if the fraction of packages that the vehicle $V_k$ contains while leaving $B_{l_r}(V_k)$ $\forall r$ is less than $w_k$ then the capacity constraint on $V_k$ is satisfied. We consider the auxiliary cost function $\sum_{k,r}\exp\big\{\theta(C(B_{l_r}(V_k),V_k)-w_k)\big\}$ where $\theta\gg 1$. Note that only when all the constraints in (\ref{eq: constLMDP}) are satisfied the auxiliary cost function value is small ($\ll 1$) otherwise it attains a large value. Similar to our approach in FLP and FLPO (Section XX and XX) we add the auxiliary cost function to the Lagrangian for the optimization problem in (\ref{eq: eqn2}) as an equality constraint that requires it to attain a specific value of $\mu_M$.  We then design our algorithm that gradually decreases $\mu_M$ till it reaches a value where all the inequality constraints (\ref{eq: constLMDP}) are satisfied and at the same time decreases the cost function $D$ corresponding to the original optimization problem. In the MEP framework, we seek to minimize the Lagrangian 
\begin{small}
\begin{align}\label{eq: LagLMDP}
\bar{F}_M &= \beta (D_M-D_{M0})\nonumber \\
&+ \beta'\Big(\sum_{k,r}\exp\big\{\theta(C(B_{l_r}(V_k),V_k)-w_k)\big\}-\mu_M\Big)- H_M,
\end{align}
\end{small}
\hspace{-0.2cm}where $D_M$, $D_{M0}$ $H_M$ are as defined in (\ref{eq: costfuncLMD}), $\beta$ and $\beta'$ are the Lagrange parameters. 
We then minimize (local) the Lagrangian $\bar{F}_M$ by setting $\frac{\partial \bar{F}_M}{\partial p_{h}(\gamma_{h+1}|\gamma_h)}=0$ to obtain the Gibbs distribution
\begin{align}
p_h(\gamma_{h+1}|\gamma_h) = e^{-\bar{c}_h}\frac{\smashoperator[r]{\sum_{\substack{(\gamma_{h+2},\hdots,\gamma_M)}}}e^{\sum_{t=h+1}^M-\bar{c}_t}}{{\sum_{\substack{(\gamma_{h+1},\hdots,\gamma_M)}}e^{\sum_{t=h}^M-\bar{c}_t}}},
\end{align}
where $\bar{c}_h := \beta c(\gamma_h,\gamma_{h+1}) + \theta\beta'e^{\theta[C(\gamma_{h+1})-w_k]}$. Since the Lagrangian $\bar{F}_M$ in (\ref{eq: LagLMDP}) is analogous to the Lagrangian $\bar{F}$ and $\bar{F}_p$ in (\ref{eq: LagFLP}) and (\ref{eq: LagFLPO}), respectively, we design annealing based algorithm similar to the Algorithm \ref{Alg: Algorithm1} and Algorithm \ref{Alg: Algorithm2} where we gradually increase $\beta$ from a small to a large value and for each $\beta$ we further anneal $\beta'$ from a small to a large value till the solution converges in the feasible region.

\section{Simulation}
\begin{figure}
    \centering
    \includegraphics[scale=0.43]{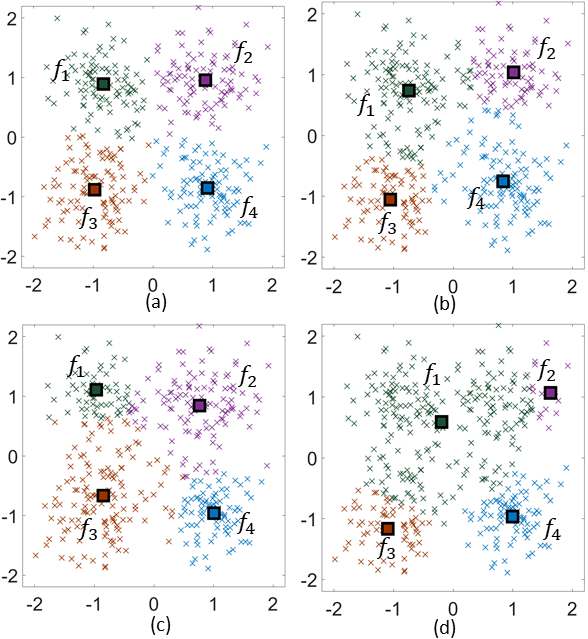}
    \caption{Illustrates unconstrained and constrained FLP (a) Unconstrained scenario, $p_1=0.24$, $p_2=0.26$, $p_3=0.26$, $p_4=0.24$ (b) $c_1=0.4$, $c_2=0.2$, $c_3=0.2$, $c_4=0.4$. Facilities allocated $p_1=0.32$, $p_2=0.19$, $p_3=0.19$, $p_4=0.30$ (c) $c_1=0.15$, $c_2=0.4$, $c_3=0.5$, $c_4=0.21$. Facilities allocated $p_1=0.15$, $p_2=0.30$, $p_3=0.34$, $p_4=0.21$ (d) $c_1=1$, $c_2=0.05$, $c_3=0.2$, $c_4=0.3$. Facilities allocated $p_1=0.59$ $p_2=0.04$, $p_3=0.16$, $p_4=0.22$ }
    \label{fig:Pic2}
\vspace{-0.7cm}
\end{figure}

In this section we illustrate our proposed framework to address the inequality based constraints in FLP, FLPO and LMDP. Figure \ref{fig:Pic2}(a) illustrates the unconstrained facility location problem with $M=4$ facilities and $N=400$ nodes randomly distributed in a $16 \text{ sq. unit}$ area. In the figure, the nodes and facilities are indicated by $\times$ and squares respectively, and the nodes are represented in the same color as the facility associated to them; for instance the green colored facility $f_1$ in Figure \ref{fig:Pic2}(a) caters to the identically colored nodes only. The final allocation of facilities is such that $p_1=0.24$, $p_2=0.26$, $p_3=0.26$, and $p_4=0.24$. In the Figure \ref{fig:Pic2}(b), the facilities $\{f_j\}$ have an inherent constraint on the capacity given as $p_j\leq c_j~\forall~j$ where $p_j$ is the effective number of nodes associated to $f_j$, $c_1=0.4$, $c_2=0.2$, $c_3=0.2$, and $c_4=0.4$. The facility allocation as given by Algorithm \ref{Alg: Algorithm1} is such that $p_1=0.30$, $p_2=0.19$, $p_3=0.19$ and $p_4=0.32$, and all the capacity constraints are satisfied. Observe that the owing to the capacity constraints, the cluster sizes for the facilities $f_2$ and $f_3$ decreases whereas it increases for the facilities $f_1$ and $f_4$. The Figures \ref{fig:Pic2}(c) and \ref{fig:Pic2}(d) consider more stringent capacity constraints on the facilities (please see Figure \ref{fig:Pic2} caption for details). For instance, in Figure \ref{fig:Pic2}(d) the facility $f_2$ requires $p_2\leq c_2=0.05$ which results into a drastic reduction in the cluster size corresponding to the facility $f_2$. Please see Figure \ref{fig:Pic2}.

\begin{figure}
    \centering    \includegraphics[scale=0.42]{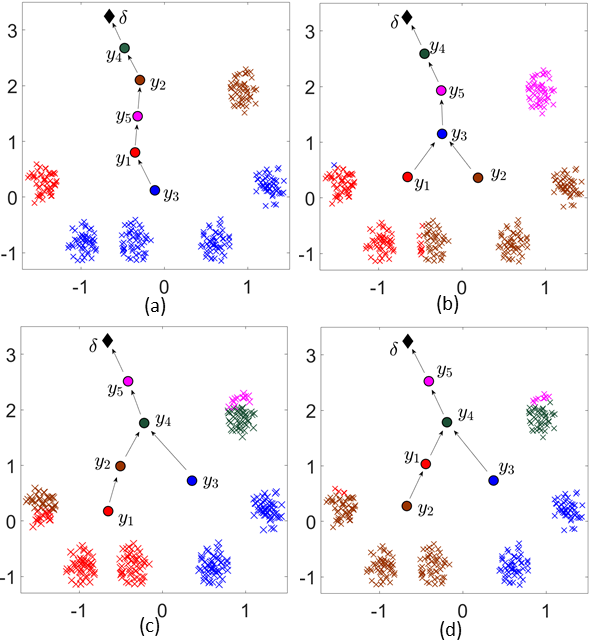}
    \caption{Illustrates unconstrained,constrained FLPO. (a) Unconstrained scenario. Facilities allocated, paths designed with $C(f_1)=0.83$, $C(f_2)=1$, $C(f_3)=0.66$, $C(f_4)=1$, $C(f_5)=0.83$, (b) Constrained scenario $w_1=0.4$, $w_2=0.8$, $w_3=w_4=w_5=1$. Facilities allocated, path designed with $C(f_1)=0.36$, $C(f_2)=0.46$, $C(f_3)=0.83$, $C(f_4)=1$, $C(f_5)=1$, (c) Constrained scenario $w_1=0.4$, $w_2=0.8$, $w_3=0.8$, $w_4=w_5=1$. Facilities allocated, path designed with $C(f_1)=0.40$, $C(f_2)=0.50$, $C(f_3)=0.33$, $C(f_4)=0.95$, $C(f_5)=1$ (d) Constrained scenario $w_1=0.5$, $w_2=0.5$, $w_3=0.8$, $w_4=w_5=1$. Facilities allocated,path designed with $C(f_1)=0.50$, $C(f_2)=0.48$, $C(f_3)=0.33$, $C(f_4)=0.97$, $C(f_5)=1$.}
    \label{fig:Pic1}
\vspace{-0.6cm}
\end{figure}

Figure \ref{fig:Pic1}(a) illustrates the unconstrained FLPO problem with $M=5$ facilities and $N=317$ nodes randomly distributed in a $12$ sq. unit area. The destination $\delta$ is denoted by a black colored diamond, the facilities are denoted by different color-filled circles and the nodes are denoted by $\times$ using a color coding that fixes the color of each node to be same as that of the facility that is the first facility on its path to the final destination $\delta$. The arrows indicate the path from each facility that leads up to the destination $\delta$. In the unconstrained setting we obtain $C(f_1)=0.83$, $C(f_2)=1$, $C(f_3)=0.66$, $C(f_4)=1$, and $C(f_5)=0.83$ where $C(f_j)$ measures the effective usage of the facility $f_j$ by the nodes $\{n_i\}$. Figures \ref{fig:Pic1}(b), (c) and (d) illustrate various instances of capacity constraints on the facilities. In the Figure \ref{fig:Pic1}(b) the capacity constraints (\ref{eq: FLPO_const}) are such that for facilities $f_1$ and $f_2$, $w_1=0.4$ and $w_2=0.8$ respectively and for the facilities $f_3,f_4$ and $f_5$ there are no capacity constraints. Since the maximum value of $C(f_j)$ is $1$ in the unconstrained setting, we set $w_3=w_4=w_5=1$ in our implementation of the algorithm. The Algorithm \ref{Alg: Algorithm2} results into facility allocation and path design such that all the capacity constraints are satisfied, i.e. $C(f_1)=0.36$, $C(f_2)=0.46$, $C(f_3)=0.83$, $C(f_4)=1$ $C(f_5)=1$ each of which is less than or equal to its upper bound $w_j$. Observe that the final facility locations and path designs by the Algorithm \ref{Alg: Algorithm2} in Figure \ref{fig:Pic1}(b) are different from the ones in Figure \ref{fig:Pic1}(a) owing to the additional inequality constraints posed by the former. For instance, in Figure \ref{fig:Pic1}(a) all the blue colored nodes (approximately $2/3^{\text{rd}}$ of the total nodes) have a corresponding $5$-hop path $y_3\rightarrow y_1\rightarrow y_5\rightarrow y_2\rightarrow y_4\rightarrow y_4\rightarrow \delta$ and the brown colored nodes follow the $2$-hop path $y_2\rightarrow y_4\rightarrow \delta$; whereas in Figure \ref{fig:Pic1}(b) all the red and brown colored nodes (approximately $5/6^{\text{th}}$ of the total nodes) follow the $4$-hop path (either $y_1\rightarrow y_3\rightarrow y_5\rightarrow y_4\rightarrow \delta$ or $y_2\rightarrow y_3\rightarrow y_5\rightarrow y_4\rightarrow \delta$) and the remaining pink colored nodes follow the $2$-hop path $y_5\rightarrow y_4\rightarrow \delta$. Please see Figure \ref{fig:Pic1} for details.

We now simulate the vehicle capacity based constraints in the Last Mile Delivery problem through a relatively small example involving four depots $B=\{B_1,B_2,B_3,B_4\}$, three service vehicles $V=\{V_1,V_2,V_3\}$ and three packages $\mathcal{P}=\{b_1,b_2,b_3\}$. The package $b_1$ originates at depot $B_1$, $b_2$ originates at depot $B_2$ and $b_3$ originates at the depot $B_3$ and all the packages are destined for the location $B_4$. The service vehicle route information is given as below
\begin{itemize}
\item $V_1$ : route is $B_1\rightarrow B_2\rightarrow B_3\rightarrow B_4$ and the respective leaving times (in minutes) are $0,30,60$, and $90$.
\item $V_2$ : route is $B_3\rightarrow B_1$ and the respective leaving time (in minutes) is $0$ and $30$.
\item $V_3$ : route is $B_2\rightarrow B_3\rightarrow B_1\rightarrow B_4$ and the respective leaving time (in minutes) is $0,20,40$, and $60$. 
\end{itemize}
In the unconstrained scenario the vehicle routes are obtained by setting $\beta\rightarrow\infty$ in the expressions of association probability (\ref{eq: eqn1}). The final routes for each of the package is given as
\begin{itemize}
\item $b_1$ : $B_1\rightarrow^{V_3}B_4$ i.e. the package $b_1$ is carried from the origin $B_1$ to the destination $B_4$ via the vehicle $V_3$. Time incurred is $60$ minutes.
\item $b_2$ : $B_2\rightarrow^{V_3}B_3\rightarrow^{V_3}B_1\rightarrow^{V_3}B_4$. Time incurred is $60$ minutes.
\item $b_3$ : $B_3\rightarrow^{V_2}B_1\rightarrow^{V_3}B_4$. Time incurred is $60$ minutes.
\end{itemize}
Note that in the above unconstrained LMDP solution, the vehicle $V_3$ carries all the three packages $b_1,b_2$ and $b_3$ in the final stretch $B_1\rightarrow B_4$. Next, we simulate the constrained LMDP case where the vehicle capacity is restricted to carry a maximum of two packages (or equivalently in our example scenario a service vehicle is allowed to carry a maximum of $2/3^{\text{rd}}$ of all the packages). In this case the our algorithm results into the same routes for the packages $b_2$ and $b_3$ as in the above unconstrained case, while the package $b_1$ is now carried from its origin $B_1$ directly to $B_4$ via the vehicle $V_1$; this ensures that the capacity constraints on the service vehicles are satisfied (vehicle $V_3$ carries only a maximum of two packages at any instant as opposed to it carrying all three packages from $B_1$ to $B_4$ in the unconstrained case); however, the time taken by $b_1$ to reach its final destination $B_4$ increases by $30$ minutes. We are working to employ our algorithm for constrained LMDP to a standard large dataset as a part of our ongoing work.

\section{CONCLUSIONS}

In this paper we propose a novel MEP based framework to incorporate inequality constraints in FLP, FLPO and LMDP. In our approach we comprehend the inequality constraints as an auxiliary cost function and use MEP to determine the associated decision variables such that the original cost function is minimized and the auxiliary cost function attains a value below a specific value corresponding to which all the inequality constraints are satisfied. The underlying idea is to design cooling laws that allows to violate the constraints and encourages exploration of the solution space  during the early stages of the algorithm; and then gradually lower the violation until the algorithm converges to a feasible point in the solution space thereby avoiding getting stuck in a poor local minima. Even though we expound on the specific optimization problems in FLP, FLPO, and LMDP, we believe our approach builds up a common framework to address constraints in many combinatorial optimization problems. 

\addtolength{\textheight}{-12cm}
\section*{APPENDIX}
\begin{itemize}
\item[1.] $A = \sum_{i=1}^M A_i$, where $A_i \in \mathbb{R}^{M\times M}$ is a diagonal matrix such that 
$(A_i)_{jj} = \sum_{\gamma_0,\gamma_1,\hdots,\gamma_{i-1}} \rho_{\gamma_0} p_0(\gamma_1|\gamma_0)\hdots p_{i-1}(f_j|\gamma_{i-1})$
\item[2.] $B = \sum_{i=1}^{M-1} (B_i + B_i^T)$ where $B_i \in \mathbb{R}^{M\times M}$ is such that 
$(B_i)_{mn} = \smashoperator{\sum\limits_{\gamma_0,\gamma_1,\hdots,\gamma_{i-1}}} \rho_{\gamma_0}p_0(\gamma_1|\gamma_0)\hdots p_{i-1}(f_m|\gamma_{i-1})p_i(f_n|f_m)$
\item[3.] $\bar{X}\in \mathbb{R}^{M\times d}$, where $\bar{X}_{mn} = \sum_{\gamma_0} \rho_{\gamma_0}p_0(f_m|\gamma_0) (\bar{z}(\gamma_0))_{n}$, where $\bar{z}(\gamma_0)_n$ is the n-th component of the spatial coordinate of $\gamma_0$.
\item[4.] $\small C = \bar{B} + \sum_{i=2}^{M-1}\underbar{B}_i + \sum_{j=2}^{M-1}\check{B}_j + \bar{D} \in \mathbb{R}^{M\times d}$, where $\small(\bar{B})_{mn} = \sum_{\gamma_0}\rho_{\gamma_0} p_0(f_m|\gamma_0)p_1(\delta|f_m)z_n$, $(\underbar{B}_i)_{mn}= \smashoperator{\sum\limits_{\gamma_0,\gamma_1,\hdots,\gamma_{i-1}}} \rho_{\gamma_0}p_0\hdots 
p_{i-2}p_{i-1}(f_m|\gamma_{i-1})p_i(\delta|f_m)z_n$ , $(\check{B}_i)_{mn}= \smashoperator{\sum_{\gamma_0,\gamma_1,\hdots,\gamma_{i-1}}} \rho_{\gamma_0} p_0\hdots p_{i-2}p_{i-1}(\delta|\gamma_{i-1})p_i(f_m|\delta)z_n$ , $(\bar{D})_{mn} = \sum_{\gamma_0,\hdots,\gamma_{M-1}}\rho_{\gamma_0}p_0(\gamma_1|\gamma_0)\hdots p_{M-1}(f_m|\gamma_{M-1})z_n$
\end{itemize}
$p_i:=p_i(\gamma_{i+1}|\gamma_i)$, $z_n$ is the $n-$th coordinate of the $\delta$.
\bibliographystyle{IEEEtran}
\bibliography{IEEEabrv}
\end{document}